\documentclass[leqno]{article}

\newtheorem{theorem}{Theorem}
\newtheorem{lemma}[theorem]{Lemma}
\newtheorem{proposition}[theorem]{Proposition}
\newtheorem{definition}[theorem]{Definition}
\newtheorem{corollary}[theorem]{Corollary}

\newcommand{\begintheorem}{\addtocounter{equation}{1}\begin{theorem}}
\newcommand{\beginlemma}{\addtocounter{equation}{1}\begin{lemma}}
\newcommand{\beginproposition}{\addtocounter{equation}{1}\begin{proposition}}
\newcommand{\begindefinition}{\addtocounter{equation}{1}\begin{definition}}
\newcommand{\begincorollary}{\addtocounter{equation}{1}\begin{corollary}}

\begin{document}

\title{Some topics in complex and harmonic analysis, 4}

\author{Stephen William Semmes	\\
	Rice University		\\
	Houston, Texas}

\date{}

\maketitle

	Let $\mathcal{S}({\bf R}^n)$ denote the Schwartz class of
rapidly decreasing smooth complex-valued functions on ${\bf R}^n$, and
let $\mathcal{S}'({\bf R}^n)$ denote the vector space of tempered
distributions on ${\bf R}^n$, which are the continuous linear
functionals on $\mathcal{S}({\bf R}^n)$.  Let $\phi \in
\mathcal{S}({\bf R}^n)$ and $\lambda \in \mathcal{S}'({\bf R}^n)$ be
given, and for each $x \in {\bf R}^n$ let $\phi_x$ be the function
defined by
\begin{equation}
	\phi_x(y) = \phi(x - y).
\end{equation}
Thus $\phi_x \in \mathcal{S}({\bf R}^n)$ for all $x \in {\bf R}^n$.

	The convolution of $\lambda$ and $\phi$ is the function
on ${\bf R}^n$ defined by
\begin{equation}
	(\lambda * \phi)(x) = \lambda(\phi_x).
\end{equation}
One can check that this is in fact a continuous function on ${\bf
R}^n$.  More precisely, $x \mapsto \phi_x$ is continuous as a function
defined on ${\bf R}^n$ with values in the Schwartz class, and hence
$\lambda(\phi_x)$ is continuous as a complex-valued function on ${\bf
R}^n$ because $\lambda$ is a continuous linear functional on the
Schwartz class.  One can show further that $\lambda * \phi$ is a
smooth function on ${\bf R}^n$.

	For each pair of multi-indices $\alpha$, $\beta$ we have the
semi-norm $\|\cdot \|_{\alpha, \beta}$ defined on the Schwartz class
by
\begin{equation}
	\|\psi\|_{\alpha, \beta} =
  \sup \biggl\{\biggl| x^\alpha \, \frac{\partial^{d(\beta)}}{\partial x^\beta}
		\, \psi(x) \biggr| : x \in {\bf R}^n \biggr\},
\end{equation}
where $d(\beta)$ denotes the degree of $\beta = (\beta_1, \ldots,
\beta_n)$, which is defined to be the sum $\beta_1 + \cdots +
\beta_n$.  This family of seminorms determines the topology on the
Schwartz class.  In particular, the statement that $\lambda$ is a
continuous linear functional on the Schwartz class means that there is
a nonnegative real number $C$ and multi-indices $\alpha_1, \beta_1,
\ldots, \alpha_p, \beta_p$ such that
\begin{equation}
	|\lambda(\psi)| \le C \sum_{j=1}^p \|\psi\|_{\alpha_j, \beta_j}
\end{equation}
for all functions $\psi$ in the Schwartz class.

	As a result, one can show that $\lambda * \phi$ is a
function of polynomial growth on ${\bf R}^n$.  To be more precise,
this means that there is a nonnegative real number $C_1$ and a
nonnegative integer $l$ such that
\begin{equation}
	|(\lambda * \phi)(x)| \le C_1 \, (1 + |x|)^l
\end{equation}
for all $x \in {\bf R}^n$.  That is, $|\lambda(\phi_x)|$ is bounded by
a constant times a finite sum of seminorms applied to $\phi_x$, and
one can check that these seminorms have polynomial growth in $x$.  One
can take this further and show that the derivatives of $\lambda *
\phi$ also have polynomial growth.

	It is perhaps better to say that $\lambda * \phi$ has at most
polynomial growth.  If $\lambda(\psi)$ is defined by integrating $\psi$
times a function in the Schwartz class, for instance, then $\lambda * \phi$
is the same as a convolution of two functions in the Schwartz class,
and is an element of the Schwartz class.  For that matter, one can show
that $\lambda * \phi$ lies in the Schwartz class if $\lambda$ has compact
support, or if $\lambda$ is defined by integration by a measure with
finite total mass.

	As another interesting special case, suppose that
$\lambda(\psi)$ is defined by integrating $\psi$ times a fixed
polynomial $p(x)$ on ${\bf R}^n$.  The definition of the convolution
$(\lambda * \phi)(x)$ is the same as the integral of $\phi(x - y)$
times $p(y)$ with respect to $y$, and this is equal to the integral of
$\phi(y)$ times $p(x - y)$ with respect to $y$.  One can expand $p(x -
y)$ as a linear combination of products of monomials in $x$ and $y$,
and the monomials in $x$ can be taken outside the integral.  Thus
$\lambda * \phi$ is a polynomial, with the coefficients given by
moments of $\phi$, which is to say integrals of $\phi$ times
monomials.  It is a simple matter to write down examples with
nonvanishing coefficients of whatever degree one might like.

	By definition, a polynomial is a linear combination of
monomials $x^\alpha$.  A polynomial of degree at most $l$ is a linear
combination of monomials $x^\alpha$ where the multi-indices $\alpha$
have degree at most $l$.  A homogeneous polynomial of degree $l$ is a
linear combination of monomials $x^\alpha$ where the multi-indices
$\alpha$ have degree equal to $l$.  This is equivalent to saying that
the polynomial $p(x)$ is homogeneous of degree $l$ in the sense that
$p(t \, x) = t^l \, p(x)$ for all $x \in {\bf R}^n$ and all positive
real numbers $t$, or arbitrary real numbers $t$ for that matter.  By
differentiating $p(t \, x)$ in $t$ and evaluating at $t = 1$ we get
the well-known identity that $\sum_{j=1}^n x_j \, (\partial / \partial
x_j) p(x)$ is equal to $l$ times $p(x)$ when $p(x)$ is homogeneous of
degree $l$.

	A polynomial $p(x)$ is said to be harmonic if $\Delta p = 0$,
where $\Delta$ denotes the usual Laplace operator, which is to say the
sum of $\partial^2 / \partial x_j^2$ over $j = 1, \ldots, n$.  If $p$
is a polynomial of degree at most $l$, then $p$ can be written in a
canonical way as the sum of polynomials which are homogeneous of
degrees ranging from $0$ to $l$.  For if $p$ is given as a linear
combination of monomials $x^\alpha$ where the multi-indices $\alpha$
have degree less than or equal to $l$, one can simply group the terms
together of the same degree.  If $p$ is a harmonic polynomial, then
the homogeneous components of $p$ are also harmonic.  For any
polynomial $p$, the Laplacian applied to the homogeneous components of
$p$ of degree $0$ or $1$ is equal to $0$, and the Laplacian applied to
a homogeneous component of $p$ of degree $m \ge 2$ is the same as the
homogeneous component of the polynomial $\Delta p$ of degree $m - 2$.

	Let $p$ be a harmonic polynomial on ${\bf R}^n$.  We would
like to check that the average of $p$ over the unit sphere in ${\bf
R}^n$ is equal to the value of $p$ at the origin.  It is enough to
check this when $p$ is a harmonic polynomial which is homogeneous of
some degree $l$, since we can express any harmonic polynomial as a sum
of homogeneous harmonic polynomials.  When $l = 0$, our homogeneous
harmonic polynomial is simply a constant, and thus its average over
the unit sphere is equal to its value at the origin.  Now suppose
that $l \ge 1$, in which case we would like to show that the average
of $p$ over the unit sphere is equal to $0$.

	By the identity mentioned earlier, if $x$ is a unit vector in
${\bf R}^n$, then $p(x)$ is equal to a nonzero constant times the
derivative of $p$ at $x$ in the direction of the unit normal to the
unit sphere.  The divergence theorem tells us that the integral of the
unit normal derivative of $p$ on the unit sphere is equal to $0$,
because the integral of the Laplacian of $p$ on the unit ball is equal
to $0$, since $p$ is harmonic.  Therefore the integral of $p$ on the
unit sphere is equal to $0$, as desired.

	This shows that the average of a harmonic polynomial on the
unit sphere is equal to the value of the polynomial at the origin.  In
general the average of a harmonic polynomial over any sphere is equal
to the value of the harmonic polynomial at the center of the sphere.
Using this one can verify that the convolution of a harmonic
polynomial $p$ with a radial function $\phi$ with suitable
integrability properties, such as a radial function in the Schwartz
class, is equal to the integral of $\phi$ over ${\bf R}^n$ times $p$.
In other words, harmonic polynomials are eigenfunctions for suitable
radial convolution operators.

\end{document}